\documentclass[a4paper,12pt]{article}

\usepackage[body={160mm,230mm}]{geometry}
\usepackage{a1}
\usepackage[english]{babel}
\usepackage[latin1]{inputenc} %permite o uso de acentos
\usepackage{graphicx}
\usepackage[dvips]{color}
\usepackage{amsfonts,amssymb}

\def\IZ{\mathbb  Z}

\def\I{\mathcal I}

\def\R{\mathcal R}

\def\mapright#1#2{\smash{\mathop{\hbox to
0.10cm{\rightarrowfill}}\limits^{#1}_{#2}}}
\def\mapleft#1#2{\smash{\mathop{\hbox to
0.10cm{\leftarrowfill}}\limits^{#1}_{#2}}}

\title{A state sum invariant for regular isotopy of links having
 	  a polynomial number of states}
\author{by Sóstenes Lins\\ {\small DMAT/CCEN/UFPE and} \\ {\small Brazilian Academy of Sciences} \\  {\small sostenes.lins@gmail.com}}
\date{Version 1, \today}

\begin{document}

\maketitle

\begin{abstract}
The state sum regular isotopy invariant of links which I introduce in this work is a generalization of the Jones Polynomial. So it distinguishes any pair of links which are distinguishable by Jones'. This new invariant, denoted {\em VSE-invariant} is strictly stronger than Jones': I detected a pair of links which are not distinguished by Jones' but are distinguished by the new invariant. The full VSE-invariant has $3^n$ states. However, there are useful specializations of it parametrized by an integer k, having $O(n^k)=\sum_{\ell=0}^k {n \choose \ell} \ 2^\ell $ states. The link with more crossings of the pair which was distinguished  by the VSE-invariant has 20 crossings. The specialization which is enough to distinguish corresponds to k=2 and has only $801$ states,  as opposed to the $2^{20} = 1,048,576$ states of the Jones polynomial of the same link. The full VSE-invariant of it has $3^{20} = 3,486,784,401$ states. The VSE-invariant is a good alternative for the Jones polynomial when the number of crossings makes the computation of this polynomial impossible. For instance, for $k=2$ the specialization of the VSE-invariant of a link with $n=500$ crossings can be computed in a few minutes, since it has only $2\,n^2+1 = 500,001$ states.
\end{abstract}

\section{Introduction: the VSE-Expansion}
The Jones polynomial, \cite{JO} or its equivalent counterpart, Kauffman's bracket \cite{KA87} does a superb job of distinguishing inequivalent knots and links. However, its computation is limited to links with a few crossing because there are $2^n$ states to be enumerated and evaluated for a link presentation having n crossings. Here I present a practical strategy to overcome exponentiability, thereby obtaining useful regular isotopy invariants with only a polynomial number of states.

\begin{figure}[h]
\begin{center}
     \includegraphics[width=12cm]{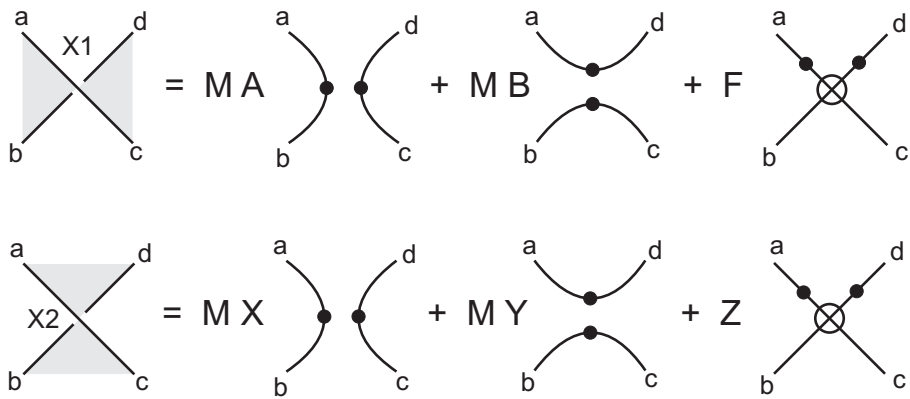}
         \caption{\bf The virtual shaded $3$-fold expansion}
     \label{fig:VSE_Expansion}
\end{center}
\end{figure}

The strategy which works is a 4-step strengthening of Kauffman's expansion for the bracket \cite{KA87}. The state sum of the VSE-invariant lives in the ring $$\R=\IZ[A, B, F, X, Y,  Z, M, o].$$
The first generalization relative to Kauffman's bracket is to use the 2-coloration (shaded and white faces) of the link diagram. This permits the distinctions of two kinds of crossing $X1$ and $X2$: the crossing of type $X1$ is the one that going counterclockwise from an overpass to an underpass the sweeped region is shaded; otherwise, if this region is white, the crossing is of type $X2$. The two types of crossings enable the definition of $4$ variables $A, B, X, Y)$ instead of the usual $2$ variables $A, B$, of only 2 of the bracket. A second generalization is that the virtual term of the expansion is included, by means of new variables $F$ and $Z$. A third generalization is to introduce a new variable $M$ to control the level: to obtain the $k$-specialization, this variable is declared to satisfy $M^{k+1}=0$. Crossings of both types are expanded according to the two rules of Fig. \ref{fig:VSE_Expansion}.

\begin{figure}[h]
\begin{center}
     \includegraphics[width=9cm]{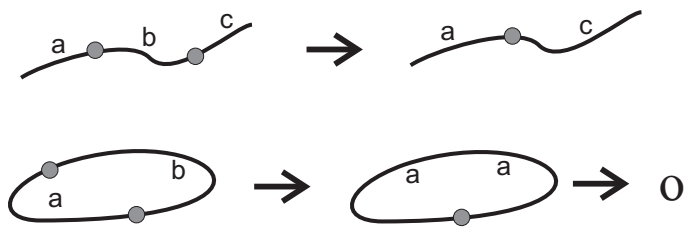}
         \caption{\bf Eliminating bivalent vertices}
     \label{fig:Simplifications}
\end{center}
\end{figure}

Note that the bracket expansion corresponds to the particular case $M=1$, $X=A$, $Y=B$ and $F=Z=0$. In the expansion some graphical bivalent vertices are created. Each monomial of a full expansion is the coefficient of a set of $m$ polygons which are evaluated by removing the bivalent vertices one by one according to Fig. \ref{fig:Simplifications}. The term $o^m$ replaces the $m$ polygons.  Variable $o$ is the {\em loop value}.

\noindent
In Mathematica, the expansion rules are given by rule1
\begin{verbatim}
rule1 = {
 X1[a_, b_, c_, d_] :> M A con[a b] con[c d] + M B con[a d] con[b c]
                                             + F con[a c] con[b d],
 X2[a_, b_, c_, d_] :> M X con[a b] con[c d] + M Y con[a d] con[b c]
                                             + Z con[a c] con[b d]
};
\end{verbatim}
The set of three simplifications, eliminating bivalent vertices is
\begin{verbatim}
rule2 = {
   con[a_ b_] con[b_ c_] :> con[a c],
   con[a_ b_] con[b_ a_] :> con[a a],
   con[a_ a_] :> o
   };
\end{verbatim}
Finally, the state sum of a product is simply
\begin{verbatim}
StateSum[product_] :=
              Expand[Simplify[(product /. rule1 // Expand) //. rule2]]
\end{verbatim}

\section{Invariance under Reidemeister moves  2 and 3}

Reidemeister move 2 can be dealt with by defining
\begin{verbatim}
LeftMove21 = X1[a, b, f, e] X2[d, e, f, c];
RightMove21 = con[a d] con[b c];
LeftMove22 = X2[a, b, f, e] X1[d, e, f, c];
RightMove22 = con[a d] con[b c];
\end{verbatim}
\begin{figure}[h]
\begin{center}
     \includegraphics[width=12cm]{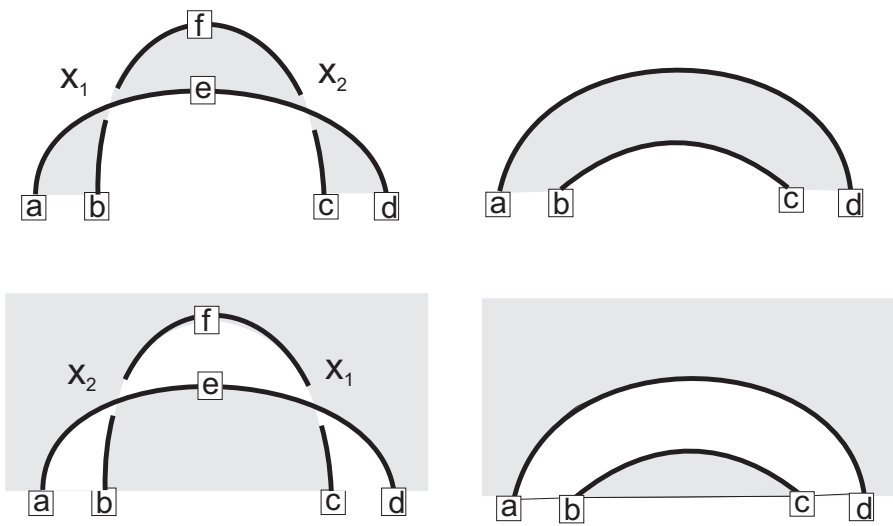}
         \caption{\bf The two types of moves 2}
     \label{fig:moves2GrayWhite}
\end{center}
\end{figure}

The fourth strengthening relative to Kauffman's bracket is consider each type of exterior to be a variable $V_{ext}$, where $ext$ is an encoding of the particular transitions relative to each type of exterior. Instead of simply imposing
$$LeftMove21-RightMove21=0, \hspace{5mm}  LeftMove22-RightMove22=0,$$
I define
$$\Omega  =  V_{con[a d] con b c]} con[a d]con[b c] + V_{con[a b] con [c d]} con[a d]con[b c] + V_{con[a c] con [cb d]} con[a c]con[b d]$$
and impose
$$ (LeftMove21-RightMove21) \ \Omega  = 0 , \hspace{5mm} (LeftMove22-RightMove22)  \ \Omega = 0.$$
\noindent
Equality must hold for all values of the exterior variables $V_{con[\_ \ \_] con [\_ \  \_]}$.
In the state sum, each exterior variable has degree at most 1.
So, if I take the partial derivatives of the state sum relative to each of these variables, the exterior variables disappear. I must impose that each such derivative must be zero, thus obtaining a polynomial equation for each exterior variables and each move. This scheme using exterior variables is clearly stronger than the usual one which does not make use of these variables: any solution of the old scheme is a solution for the new scheme but not vice-versa. The six equations coming from Reidemeister move 2 are:
\begin{center}
$eq_1:= B o^2 X M^2+A o X M^2+A o^2 Y M^2+B o Y M^2+F o X M+F o Y M+A o Z
   M+B o Z M-o^2+F o^2 Z = 0,$\\
$eq_2:= A M^2 Y o^3+A M^2 X o^2+B M^2 Y o^2+F M Y o^2+A M Z o^2+B M^2 X
   o+F M X o+F Z o+B M Z o-o=0,$\\
$eq_3:= A o X M^2+B o X M^2+A o^2 Y M^2+B o Y M^2+F o^2 X M+F o Y M+B o^2
   Z M+A o Z M-o+F o Z=0,$\\
$eq_4:= -o^2 + F M o X + A M^2 o X + B M^2 o^2 X + F M o Y + B M^2 o Y +
 A M^2 o^2 Y + A M o Z + B M o Z + F o^2 Z=0,$\\
 $eq_5:= -o + F M o^2 X + A M^2 o^2 X + B M^2 o^3 X + F M o Y + A M^2 o Y +
 B M^2 o^2 Y + F o Z + A M o Z + B M o^2 Z=0,$\\
$eq_6:=-o + F M o X + A M^2 o X + B M^2 o^2 X + A M^2 o Y + B M^2 o Y +
F M o^2 Y + F o Z + $\\
$ B M o Z + A M o^2 Z=0.$
\end{center}
Note that due to symmetry, $eq_1$ and $eq_4$ coincide and there are only 5 distinct relations.

\begin{figure}[h]
\begin{center}
     \includegraphics[width=12cm]{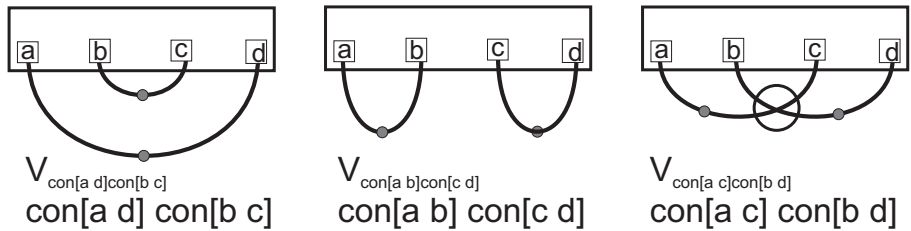}
         \caption{\bf Three types of exteriors for Reidemeister move 2}
     \label{fig:Exterior2}
\end{center}
\end{figure}
\begin{figure}[h]
\begin{center}
     \includegraphics[width=12cm]{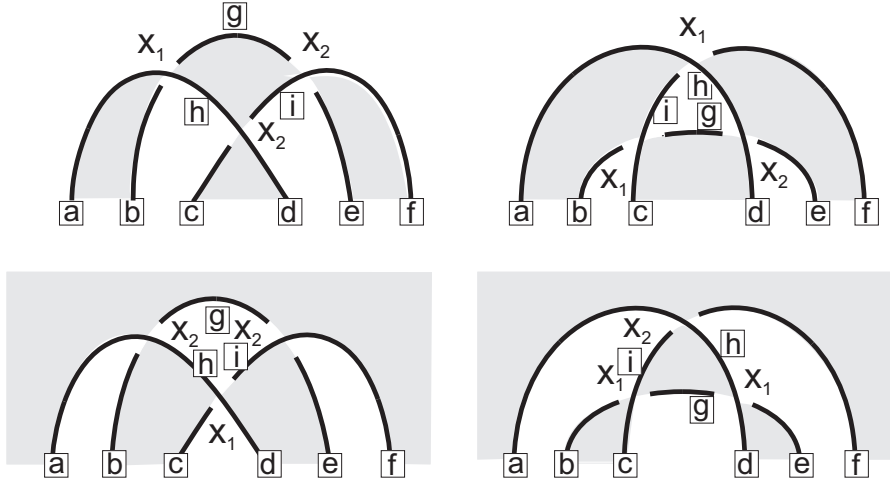}
         \caption{\bf The two types of moves 3}
     \label{fig:moves3GrayWhite}
\end{center}
\end{figure}

\noindent
Now I do a similar job for move 3. For each such move there are now 15 exterior variables, giving rise to 15 equations. For moves 2 and 3 there is a total of 2(3+15)=36 equations in the 8 variables, but there are only 27 distinct ones,  due to symmetry. The encoding in Mathematica of the two types of moves 3 are
\begin{verbatim}
LeftMove31 = X1[a, b, h, g] X2[i, e, f, g] X2[h, c, d, i];
RightMove31 = X1[a, i, h, f] X1[c, g, i, b] X2[h, g, d, e] ;
LeftMove32 = X2[a, b, h, g] X1[i, e, f, g] X1[h, c, d, i];
RightMove32 = X2[a, i, h, f] X2[c, g, i, b] X1[h, g, d, e];
\end{verbatim}
These moves correspond to the situation of Fig. \ref{fig:moves3GrayWhite}.

\section{Obtaining the relevant ideal}
Instead of considering a set of 27 polynomial equations $pol_i=0$, in the spirit of King, \cite{King2005}, I take the ideal generated by the left hand side of the system of equations.  These polynomials generate an ideal, named  $\I_\infty$. Instead of solving the system of polynomial equations, I compute a Gr\"obner for the ideal $\I_\infty$ relative to a fixed monomial ordering. The VSE-invariant is the normal form $\eta_\infty(p)$ of the classes of polynomials $\overline{p} \in \R/\I_\infty$. If $p$ and $q$ are VSE-state sums of two links $L_p$ and $L_q$ which can be transformed one into the other by Reidemeister moves 2 and 3,
then $\eta_\infty(p) = \eta_\infty(q)$.

\subsection{The ideal $\I_\infty$}
I have written a subroutine
in Mathematica to obtain automatically the polynomials relative to a given set of moves.
The ideal $\I_\infty$ of $\R=\IZ[A, B, F, X, Y,  Z, M, o]$ corresponding to Reidemeister moves 2 and 3 is
$$\I_\infty = \langle pol_1, pol_2,\ldots,pol_{26}, pol_{27} \rangle,$$
where\\
\noindent
$pol_1=o (B o^2 X M^2+B o Y M^2+B o Z M+A (M o X+M Y+Z) M+F(M (o X+Y)+Z)-1),$\\
   $pol_2= o (A o X M^2+A o^2 Y M^2+A o Z M+B(M X+M o Y+Z) M+F (M (X+o Y)+Z)-1),$\\
   $pol_3=o (M (A M X+B M Y+F (X+Y)+A Z+B Z)+o (B X M^2+A Y M^2+F Z-1)),$\\
   $pol_4=o (B o X M^2+B Y M^2+B Z M+A (M (X+Y)+o Z) M+F (M (X+o Y)+Z)-1),$\\
   $pol_5=o (F M o X+F M Y+F Z+A M (M (X+o Y)+Z)+B M (M (X+Y)+o Z)-1),$\\
   $pol_6=o (A^2 (M (o X+Y)+Z) M^2+B^2 (M (o X+Y)+Z) M^2-B ((X^2+2 o Y X+Y^2) M^2+2 (X+Y) Z M+o Z^2-2 F (M (X+Y)+o Z)) M-A (o X^2 M^2+o Y^2
   M^2+2 X Y M^2+2 X Z M+2 Y Z M-2 B (M X+M o Y+Z) M+o Z^2-2 F (M  (X+Y)+o Z)) M+F (F o-Z o-M (X+Y)) (M (X+Y)+o Z)),$\\
   $pol_7=o (-A^2 (M (o X+Y)+Z) M^2-B^2 (M (o X+Y)+Z) M^2+B((X^2+2 o Y X+Y^2) M^2+2 (X+Y) Z M+o Z^2-2 F (M
   (X+Y)+o Z)) M+A (o X^2 M^2+o Y^2 M^2+2 X Y M^2+2 X Z
   M+2 Y Z M-2 B (M X+M o Y+Z) M+o Z^2-2 F (M (X+Y)+o Z))
   M-F (F o-Z o-M (X+Y)) (M (X+Y)+o Z)),$\\
   $ pol_8=o (B^2 (M (o
   X+Y)+Z) M^2+A^2 o (M (X+o Y)+Z) M^2-A (M (X+o Y)+Z) (-2 F-2 B
   M+o (M X+M o Y+Z)) M-B ((M (X+o Y)+Z)^2-2 F (M (X+Y)+o
   Z)) M+F (F (M (o X+Y)+Z)-(M (X+o
   Y)+Z)^2)),$\\
   $pol_9=o (-B^2 (M (o X+Y)+Z) M^2-A^2 o (M
   (X+o Y)+Z) M^2+A (M (X+o Y)+Z) (-2 F-2 B M+o (M X+M o Y+Z)) M+B
   ((M (X+o Y)+Z)^2-2 F (M (X+Y)+o Z)) M+F ((M
   (X+o Y)+Z)^2-F (M (o X+Y)+Z))),$\\
   $ pol_{10}=o (A^2 (M (o
   X+Y)+Z) M^2+B^2 (M (o X+Y)+Z) M^2-A (o X^2 M^2+o Y^2 M^2+2
   X Y M^2+2 X Z M+2 Y Z M-2 B (M (X+Y)+o Z) M+o Z^2-2 F (M (X+o
   Y)+Z)) M-B ((M (X+Y)+o Z)^2-2 F (M (X+o Y)+Z))
   M+F (-(X^2+2 o Y X+Y^2) M^2-2 (X+Y) Z M-o Z^2+F
   o (M (X+o Y)+Z))),$\\
   $ pol_{11}=o (-A^2 (M (o X+Y)+Z)
   M^2-B^2 (M (o X+Y)+Z) M^2+A (o X^2 M^2+o Y^2 M^2+2 X Y
   M^2+2 X Z M+2 Y Z M-2 B (M (X+Y)+o Z) M+o Z^2-2 F (M (X+o
   Y)+Z)) M+B ((M (X+Y)+o Z)^2-2 F (M (X+o Y)+Z))
   M+F ((X^2+2 o Y X+Y^2) M^2+2 (X+Y) Z M+o Z^2-F
   o (M (X+o Y)+Z))),$\\
   $pol_{12}=o (B^2 (M (o X+Y)+Z) M^2+A^2
   o (M (X+Y)+o Z) M^2-B ((X^2+2 Y X+o Y^2) M^2+2
   (o X+Y) Z M+Z^2-2 F (M (X+o Y)+Z)) M-A (-2 F (M
   (X+Y)+o Z)-2 B M (M (X+Y)+o Z)+o ((X^2+2 Y X+o
   Y^2) M^2+2 (o X+Y) Z M+Z^2)) M+F
   (-(X^2+2 Y X+o Y^2) M^2-2 (o X+Y) Z M-Z^2+F (M
   (o X+Y)+Z))),$ \\
   $pol_{13}=o (-B^2 (M (o X+Y)+Z) M^2-A^2 o
   (M (X+Y)+o Z) M^2+B ((X^2+2 Y X+o Y^2) M^2+2 (o
   X+Y) Z M+Z^2-2 F (M (X+o Y)+Z)) M+A (-2 F (M (X+Y)+o
   Z)-2 B M (M (X+Y)+o Z)+o ((X^2+2 Y X+o Y^2)
   M^2+2 (o X+Y) Z M+Z^2)) M+F ((X^2+2 Y X+o
   Y^2) M^2+2 (o X+Y) Z M+Z^2-F (M (o
   X+Y)+Z))),$\\
   $pol_{14}=o (A^2 o (M (o X+Y)+Z) M^2+B^2 o (M
   (o X+Y)+Z) M^2-B ((2 X Y+o
   (X^2+Y^2)) M^2+2 (X+Y) Z M+o Z^2-2 F (M (o
   X+Y)+Z)) M-A (-2 F (M (o X+Y)+Z)-2 B M (M (o
   X+Y)+Z)+o ((2 X Y+o (X^2+Y^2)) M^2+2
   (X+Y) Z M+o Z^2)) M+F (-(2 X Y+o
   (X^2+Y^2)) M^2-2 (X+Y) Z M-o Z^2+F o (M (o
   X+Y)+Z))),$\\
   $pol_{15}=o (-A^2 o (M (o X+Y)+Z) M^2-B^2 o (M
   (o X+Y)+Z) M^2+B ((2 X Y+o
   (X^2+Y^2)) M^2+2 (X+Y) Z M+o Z^2-2 F (M (o
   X+Y)+Z)) M+A (-2 F (M (o X+Y)+Z)-2 B M (M (o
   X+Y)+Z)+o ((2 X Y+o (X^2+Y^2)) M^2+2
   (X+Y) Z M+o Z^2)) M+F ((2 X Y+o
   (X^2+Y^2)) M^2+2 (X+Y) Z M+o Z^2-F o (M (o
   X+Y)+Z))),$\\
   $pol_{16}=o (-(M (o X+Y)+Z) F^2+(o X^2
   M^2+Y^2 M^2+2 X Y M^2+2 X Z M+2 o Y Z M-2 B o (M o X+M Y+Z) M-2
   A (M (o X+Y)+Z) M+Z^2) F+M (-M (M (o X+Y)+Z)
   A^2+((2 X Y+o (X^2+Y^2)) M^2+2 (X+Y)
   Z M-2 B o (M (o X+Y)+Z) M+o Z^2) A-B (B M o^2-M X
   o-M Y-Z) (M (o X+Y)+Z))),$\\
   $pol_{17}=o (-(M (o
   X+Y)+Z) F^2-(2 B M+2 A o M-o X M-Y M-Z) (M (o X+Y)+Z) F+M
   (-M (M (o X+Y)+Z) A^2+((2 X Y+o
   (X^2+Y^2)) M^2+2 (X+Y) Z M-2 B (M (o X+Y)+Z)
   M+o Z^2) A+B ((o X^2+2 Y X+Y^2) M^2+2
   (X+o Y) Z M-B o (M (o X+Y)+Z) M+Z^2))),$\\
   $pol_{18}=o((M (o X+Y)+Z) F^2+(2 B M+2 A o M-o X M-Y M-Z) (M (o
   X+Y)+Z) F+M (M (M (o X+Y)+Z) A^2+(-(2 X Y+o
   (X^2+Y^2)) M^2-2 (X+Y) Z M+2 B (M (o X+Y)+Z)
   M-o Z^2) A+B (-(o X^2+2 Y X+Y^2) M^2-2
   (X+o Y) Z M+B o (M (o X+Y)+Z) M-Z^2))),$\\
   $ pol_{19}=o
   ((M (o X+Y)+Z) F^2+(-o X^2 M^2-Y^2 M^2-2 X Y M^2-2 X
   Z M-2 o Y Z M+2 B o (M o X+M Y+Z) M+2 A (M (o X+Y)+Z)
   M-Z^2) F+M (M (M (o X+Y)+Z) A^2+(-(2 X Y+o
   (X^2+Y^2)) M^2-2 (X+Y) Z M+2 B o (M (o X+Y)+Z)
   M-o Z^2) A+B (B M o^2-M X o-M Y-Z) (M (o
   X+Y)+Z))),$\\
   $pol_{20}=o (-(M (X+Y)+o Z) F^2+(o X^2
   M^2+Y^2 M^2+o X Y M^2+X Y M^2+o^2 X Z M+X Z M+o Y Z M+Y Z M-A
   (Z o^2+2 M X o+M Y o+M Y+Z) M-B (M (2 X+o Y+Y)+(o+1)
   Z) M+o Z^2) F+M (-M (M (X+Y)+o Z)
   A^2+((X^2+2 Y X+o Y^2) M^2+2 (o X+Y) Z M-B (M
   (2 X+o Y+Y)+(o+1) Z) M+Z^2) A+B ((X^2+(o+1) Y
   X+Y^2) M^2+(o+1) (X+Y) Z M-B (M (o X+Y)+Z)
   M+Z^2))), $\\
   $pol_{21}=o (-(M (X+o Y)+Z) F^2+(o
   X^2 M^2+o Y^2 M^2+o^2 X Y M^2+X Y M^2+o X Z M+X Z M+o Y Z M+Y Z
   M-A (M Y o^2+2 M X o+Z o+M Y+Z) M-B (M (2 X+o
   Y+Y)+(o+1) Z) M+Z^2) F+M (-M (M (X+o Y)+Z)
   A^2+((M (X+o Y)+Z)^2-B M (M (2 X+o Y+Y)+(o+1) Z))
   A+B ((X+Y) (X+o Y) M^2+(Y o^2+X o+X+Y) Z M-B (M
   (o X+Y)+Z) M+o Z^2))),$\\
   $pol_{22}=o ((M (X+Y)+o Z)
   F^2+(-o X^2 M^2-Y^2 M^2-o X Y M^2-X Y M^2-o^2 X Z M-X Z
   M-o Y Z M-Y Z M+A (Z o^2+2 M X o+M Y o+M Y+Z) M+B (M
   (2 X+o Y+Y)+(o+1) Z) M-o Z^2) F+M (M (M (X+Y)+o Z)
   A^2+(-(X^2+2 Y X+o Y^2) M^2-2 (o X+Y) Z M+B (M
   (2 X+o Y+Y)+(o+1) Z) M-Z^2) A+B (-(X^2+(o+1) Y
   X+Y^2) M^2-(o+1) (X+Y) Z M+B (M (o X+Y)+Z)
   M-Z^2))), $\\
   $pol_{23}=o ((M (X+o Y)+Z) F^2+(-o
   X^2 M^2-o Y^2 M^2-o^2 X Y M^2-X Y M^2-o X Z M-X Z M-o Y Z M-Y Z
   M+A (M Y o^2+2 M X o+Z o+M Y+Z) M+B (M (2 X+o
   Y+Y)+(o+1) Z) M-Z^2) F+M (M (M (X+o Y)+Z)
   A^2+(B M (M (2 X+o Y+Y)+(o+1) Z)-(M (X+o Y)+Z)^2)
   A+B (-(X+Y) (X+o Y) M^2-(Y o^2+X o+X+Y) Z M+B
   (M (o X+Y)+Z) M-o Z^2))), $\\
   $pol_{24}=o (-(M (X+Y)+o
   Z) F^2+(X^2 M^2+o Y^2 M^2+o X Y M^2+X Y M^2+o X Z M+X Z
   M+o^2 Y Z M+Y Z M-B (Z o^2+2 M X o+M Y o+M Y+Z) M-A
   (M (2 X+o Y+Y)+(o+1) Z) M+o Z^2) F+M (-M (M (X+o
   Y)+Z) A^2+((M (X+o Y)+Z)^2-B M (M (Y o^2+2 X
   o+Y)+(o+1) Z)) A-B (B M o-M Y o-M X-Z) (M (o
   X+Y)+Z))), $\\
   $pol_{25}=o ((M (X+Y)+o Z) F^2+(-X^2
   M^2-o Y^2 M^2-o X Y M^2-X Y M^2-o X Z M-X Z M-o^2 Y Z M-Y Z M+B
   (Z o^2+2 M X o+M Y o+M Y+Z) M+A (M (2 X+o Y+Y)+(o+1)
   Z) M-o Z^2) F+M (M (M (X+o Y)+Z) A^2+(B M
   (M (Y o^2+2 X o+Y)+(o+1) Z)-(M (X+o
   Y)+Z)^2) A+B (B M o-M Y o-M X-Z) (M (o
   X+Y)+Z))), $\\
   $pol_{26}=o (-(M (X+o Y)+Z) F^2+(X^2
   M^2+Y^2 M^2+o X Y M^2+X Y M^2+o X Z M+X Z M+o Y Z M+Y Z M-B
   (M Y o^2+2 M X o+Z o+M Y+Z) M-A (M (2 X+o Y+Y)+(o+1)
   Z) M+Z^2) F+M (-M (M (X+Y)+o Z)
   A^2+((X^2+2 Y X+o Y^2) M^2+2 (o X+Y) Z M-B
   (M (2 o X+o Y+Y)+(o^2+1) Z) M+Z^2)
   A-B (M (o X+Y)+Z) (B M o-Z o-M (X+Y)))), $\\
   $pol_{27}=o ((M
   (X+o Y)+Z) F^2+(-X^2 M^2-Y^2 M^2-o X Y M^2-X Y M^2-o X Z
   M-X Z M-o Y Z M-Y Z M+B (M Y o^2+2 M X o+Z o+M Y+Z)
   M+A (M (2 X+o Y+Y)+(o+1) Z) M-Z^2) F+M (M (M (X+Y)+o
   Z) A^2+(-(X^2+2 Y X+o Y^2) M^2-2 (o X+Y) Z M+B
   (M (2 o X+o Y+Y)+(o^2+1) Z) M-Z^2)
   A+B (M (o X+Y)+Z) (B M o-Z o-M (X+Y)))).$

\subsection{The Gr\"obner basis $B_\infty$}
The Gr\"obner basis $B_\infty$ for $\I_\infty$ relative to the lexicographical order of the monomials
with variable order $vars=\{A, B, F, X, Y, Z, M, o\}$ has only 15 polynomials: \vspace{3mm} \\
$B_\infty=\{p_1=Z^3 o^4-Z o^4-Z^3 o^3+Z o^3-4 Z^3 o^2+4 Z o^2+4 Z^3 o-4 Z
   o,$\\
   $p_2=-Z^3 o^3+M^2 Y^2 Z o^3+Z o^3-Z^3 o^2+M^2 Y^2 Z o^2+Z o^2+2
   Z^3 o-2 M^2 Y^2 Z o-2 Z o,$\\
   $p_3=M^4 o^2 Y^4-M^4 o Y^4-4 M^3 o^2 Z
   Y^3+4 M^3 o Z Y^3+M^2 o^3 Y^2-M^2 o^2 Y^2+6 M^2 o^2 Z^2 Y^2-6
   M^2 o Z^2 Y^2+2 M o^3 Z^3 Y-2 M o^2 Z^3 Y-2 M o^3 Z Y+2 M o^2 Z
   Y-2 o^3 Z^4-o^2 Z^4+3 o Z^4+o^2+2 o^3 Z^2-2 o Z^2-o,$\\
   $p_4=Z^3 o^3-M Y
   Z^2 o^3+M Y o^3-Z o^3+M^3 Y^3 o^2+2 M Y Z^2 o^2+M X o^2-M Y
   o^2-3 M^2 Y^2 Z o^2-M^3 Y^3 o-Z^3 o-M Y Z^2 o-M X o+3 M^2 Y^2 Z
   o+Z o,$\\
   $p_5=Y Z^2 o^4-Y o^4-M^2 Y^3 o^3+X Z^2 o^3+Y Z^2 o^3-X o^3-Y
   o^3-M^2 Y^3 o^2+X Z^2 o^2-2 Y Z^2 o^2-X o^2+2 Y o^2+2 M^2 Y^3
   o-2 X Z^2 o+2 X o,$\\
   $p_6=o^3 Y Z^6-o Y Z^6-M o^3 Y^2 Z^5+2 M o^2 Y^2
   Z^5-M o Y^2 Z^5-3 M^2 o^2 Y^3 Z^4+3 M^2 o Y^3 Z^4+o^2 X Z^4-o X
   Z^4-o^3 Y Z^4-o^2 Y Z^4+2 o Y Z^4+M^3 o^2 Y^4 Z^3-M^3 o Y^4
   Z^3+M o^3 Y^2 Z^3-M o^2 Y^2 Z^3+2 M^2 o^2 Y^3 Z^2-2 M^2 o Y^3
   Z^2-2 o^2 X Z^2+2 o X Z^2-o^3 Y Z^2+2 o^2 Y Z^2-o Y Z^2-M^3 o^2
   Y^4 Z+M^3 o Y^4 Z-M o^2 Y^2 Z+M o Y^2 Z+M^2 o^2 Y^3-M^2 o
   Y^3+o^2 X-o X+o^3 Y-o^2 Y,$\\
   $p_7=-2 o^3 Z^4+2 o Z^4+2 M o^3 Y Z^3-4 M
   o^2 Y Z^3+2 M o Y Z^3+2 o^3 Z^2-o^2 Z^2+6 M^2 o^2 Y^2 Z^2-6 M^2
   o Y^2 Z^2-2 M^3 o^2 Y^3 Z+2 M^3 o Y^3 Z+2 M o X Z-2 M o^3 Y Z+4
   M o^2 Y Z+o^2+M^2 o X^2+M^2 o Y^2-2 o+2 M^2 o X Y,$\\
   $p_8=F o^3-Z o^3+F
   o^2-Z o^2-2 F o+2 Z o,$\\
   $p_9=o^2 Y Z^5-o Y Z^5-3 M o^2 Y^2 Z^4+3 M o
   Y^2 Z^4-F o^2 Y Z^4+F o Y Z^4+3 M^2 o^2 Y^3 Z^3-3 M^2 o Y^3
   Z^3+3 F M o^2 Y^2 Z^3-3 F M o Y^2 Z^3-o^2 X Z^3+o X Z^3+o^2 Y
   Z^3-o Y Z^3-M^3 o^2 Y^4 Z^2+M^3 o Y^4 Z^2-3 F M^2 o^2 Y^3 Z^2+3
   F M^2 o Y^3 Z^2-M o^2 Y^2 Z^2+M o Y^2 Z^2+F o^2 X Z^2-F o X
   Z^2-F o^2 Y Z^2+F o Y Z^2+F M^3 o^2 Y^4 Z-F M^3 o Y^4 Z+M^2 o^2
   Y^3 Z-M^2 o Y^3 Z+F M o^2 Y^2 Z-F M o Y^2 Z+o^2 X Z-o X Z-2 o^2
   Y Z+2 o Y Z-F M^2 o^2 Y^3+F M^2 o Y^3-F o^2 X+F o X+2 F o^2 Y-2
   F o Y,$\\
   $p_{10}=o^2 X^2 F^2-o X^2 F^2+o^2 Y^2 F^2-o Y^2 F^2-2 o^2 X Y
   F^2+2 o X Y F^2-2 o^2 X^2 Z F+2 o X^2 Z F-2 o^2 Y^2 Z F+2 o Y^2
   Z F+4 o^2 X Y Z F-4 o X Y Z F+o^2 X^2 Z^2-o X^2 Z^2+o^2 Y^2
   Z^2-o Y^2 Z^2-2 o^2 X Y Z^2+2 o X Y Z^2,$\\
   $p_{11}=-F o^2+B M o^2-M Y
   o^2+Z o^2+F o-B M o+M Y o-Z o,$\\
   $p_{12}=-F M^3 o^2 Y^4+F M^3 o Y^4+M^3
   o^2 Z Y^4-M^3 o Z Y^4-3 M^2 o^2 Z^2 Y^3+3 M^2 o Z^2 Y^3+3 F M^2
   o^2 Z Y^3-3 F M^2 o Z Y^3+3 M o^2 Z^3 Y^2-3 M o Z^3 Y^2+2 F M
   o^2 Y^2-3 F M o^2 Z^2 Y^2+3 F M o Z^2 Y^2-2 F M o Y^2-2 M o^2 Z
   Y^2+2 M o Z Y^2-o^2 Z^4 Y+o Z^4 Y+F o^2 Z^3 Y-F o Z^3 Y+o^2
   Y+o^2 Z^2 Y-o Z^2 Y-o Y-2 F o^2 Z Y+2 F o Z Y-B o^2+B o^2 Z^2-B
   o Z^2+o^2 X Z^2-o X Z^2+B o-F o^2 X Z+F o X Z,$\\
   $p_{13}=F o^2-Z o^2+A M o+B M o-M X o-M Y o,$\\
   $p_{14}=2 F M^3 o^2 Y^4-2 F M^3 o Y^4-2 M^3 o^2 Z
   Y^4+2 M^3 o Z Y^4+6 M^2 o^2 Z^2 Y^3-6 M^2 o Z^2 Y^3-6 F M^2 o^2
   Z Y^3+6 F M^2 o Z Y^3-6 M o^2 Z^3 Y^2+6 M o Z^3 Y^2-4 F M o^2
   Y^2+6 F M o^2 Z^2 Y^2-6 F M o Z^2 Y^2+3 F M o Y^2+4 M o^2 Z
   Y^2-3 M o Z Y^2+2 o^2 Z^4 Y-2 o Z^4 Y-2 F o^2 Z^3 Y+2 F o Z^3
   Y-2 o^2 Z^2 Y+3 o Z^2 Y+o Y-2 F M o X Y+2 F o^2 Z Y-4 F o Z Y+2
   M o X Z Y-F M o X^2+A o Z^2+B o Z^2+o X Z^2-A o-B o+o X+M o X^2
   Z-2 F o X Z,$\\
   $p_{15}=A F o^2-B F o^2-F X o^2+F Y o^2-A Z o^2+B Z o^2+X Z
   o^2-Y Z o^2-A F o+B F o+F X o-F Y o+A Z o-B Z o-X Z o+Y Z
   o\}.$

This basis is obtained with the Mathematica command
$$B_\infty = GroebnerBasis[\I_\infty, vars],$$
where $\I_\infty=\{pol_1, pol_2, \ldots, pol_{26}, pol_{27}\}$. The normal form of a polynomial $poly$ relative to the Gr\"obner basis $B_\infty$ is obtained in Mathematica
simply by writing
$$\eta_\infty[poly\_] :=
  PolynomialReduce[poly, B_\infty, vars][[2]];$$
This normal form is a regular isotopy invariant of links which generalizes the Jones polynomial. The {\em VSE-invariant of a link} is defined to be the normal form $\eta_\infty$ relative to the Gr\"obner basis $GB_\infty$ applied to the $VSE$-state sum of the link.

\subsection{The $k$-specializations of the VSE-invariant}

Define $\I_k = \I_\infty \cup \{M^{k+1}\}$ and let $B_k$ be the Gr\"obner basis for $\I_k$ with the same monomial order as before. I have computed explicitly $B_1, B_2, \ldots, B_{10}, B_{11}$. They have respectively $14, 25, 30, 37, 44, 53, 62, 73, 84, 97, 110$ (not horrendous) polynomials, therefore explicit computations via the normal forms are available. The corresponding normal forms $\eta_1, \eta_2, \ldots, \eta_{10}, \eta_{11}, \ldots$ are regular isotopy invariants of a link $L$ when applied to their VSE-state sum $ss(L)$. The proof of the following proposition is straightforward:
\begin{proposition}
The number of non-null states for computing the regular isotopy invariant $\eta_k(ss(L_n))\equiv\eta_k(L_n)$, of a link $L_n$ with $n$ crossings is
$\sum_{\ell=0}^k {n \choose \ell} \ 2^\ell,$ for $k \le n$ and $3^n$ for $k>n$.
\end{proposition}

\section{Comparing the VSE-invariant with the bracket}

In this section we present examples of computations of the VSE-invariant on some knots and links.
It seems, from these examples that the $VSE$-invariants and the bracket have always the same discriminative power. However... see last subsection!
For the notation on the knots see \cite{DBN}.
\subsection{Knots $4_1$ and $K11n19$}
The knots $4_1$ and $K11n19$ with with writhes $-3$ have the same $\eta_\infty$ invariant.
This implies that $\eta_k(4_1)=\eta_k(K11n19)$ for all integer values $k>0$.
\begin{figure}[h]
\begin{center}
     \includegraphics[width=15cm]{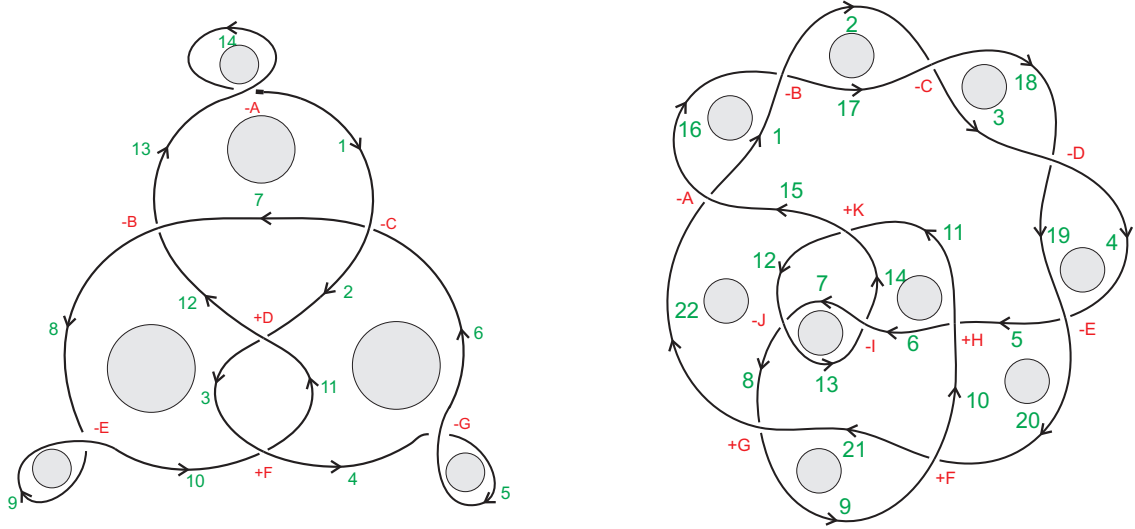}
         \caption{\bf Jones and VSE-invariant agree on $4_1$ (with write -3) and $K11n19$}
     \label{fig:4_1W3K11n19}
\end{center}
\end{figure}

\noindent $\eta_1=-\frac{1}{2} o (6 M (o-1) Y+Z (4 Z^2-3 o
   (Z^2-1)-6)),$

\noindent $\eta_2=\frac{1}{8} o (Z (12 Z^4-52 Z^2+6 o^2
   (Z^2-1)-9 o (Z^4-4 Z^2+3)+48)-24 M
   (o-1) Y),$

\noindent $\eta_3=\frac{1}{16} o (Z (12 (Z^2-1) o^2+3 (5
   Z^6-21 Z^4+39 Z^2-23) o+$\\
   .\hspace{15mm} $4 (-5 Z^6+21 Z^4-41
   Z^2+29))-16 M (o-1) Y (2 Z^2+o
   (Z^2-1)+1)),$

\noindent $\eta_4=\frac{1}{128} o (Z (24 (Z^4+2 Z^2-3) o^2-3
   (35 Z^8-180 Z^6+370 Z^4-436 Z^2+211) o+$\\
   .\hspace{15mm} $4 (35
   Z^8-180 Z^6+366 Z^4-444 Z^2+255))-128 M (o-1) Y
   (2 Z^2+o (Z^2-1)+1)),$\\
   .\\
\noindent $\eta_\infty=o (M^3 (o^8-8 o^6+21 o^4-20 o^2+6) Y^3+6 M^2 (o-1)
   Z Y^2+Z (-(Z^2-1) o^2+$\\
   .\hspace{15mm} $(Z^2-1)
   o+1)+M (X+Y (-(Z^2-1) o^9+9
   (Z^2-1) o^7-28 (Z^2-1) o^5+$\\
   .\hspace{15mm} $35
   (Z^2-1) o^3+(16-19 Z^2) o+4
   Z^2))).$

\subsection{Knots $8_8$ and $Mirror 10_{129}$}
\begin{figure}[h]
\begin{center}
     \includegraphics[width=14cm]{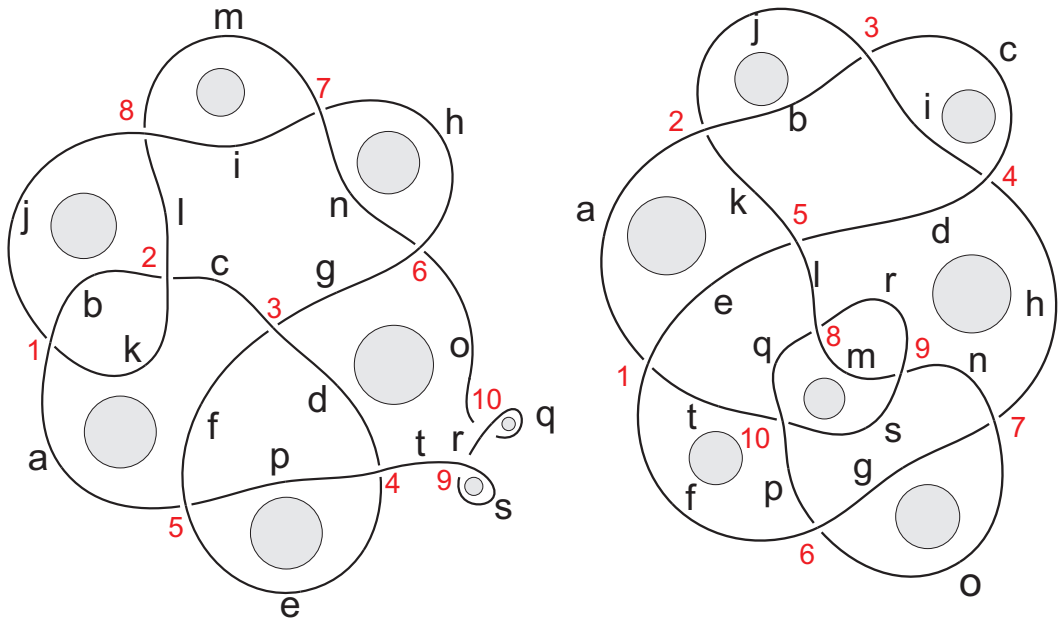}
         \caption{\bf Jones and VSE-invariant agree on $8_8$ (with writhe +2) and $Mirror 10_{129}$}
     \label{fig:8_8and10_129}
\end{center}
\end{figure}

\noindent For this pair of knots, the specializations (computed from k=1 up to 11) are rather insensitive:\\
\noindent $\eta_k(8_8)=o=\eta_k(10_{129}), k=1,2,\ldots,10,11.$ At the $\infty$-level:

\noindent $\eta_\infty(8_8)=o (M^2 Y^2 o^9-(Z^2-1) o^8-10 M^2 Y^2 o^7+10
   (Z^2-1) o^6+35 M^2 Y^2 o^5-$\\   .\hspace{15mm} $35 (Z^2-1)
   o^4-50 M^2 Y^2 o^3+50 (Z^2-1) o^2+24 M^2 Y^2 o-24
   Z^2+25)=\eta_\infty(10_{129}).$

\subsection{Knots $9_{42}$ and $Mirror 9_{42}$}
The knots $9_{42}$ and $Mirror 9_{42}$ are indistinguishable by the Jones invariant, by the Kauffman invariant and by the HOMFLY invariant. The full VSE-invariant also does not distinguishes them. As with the
previous pair, the specializations collapses to $o$ for $k$ from 1 to $11$. At $k=\infty$,
$$\eta_\infty(9_{42})=o (-7 + 8 Z^2 - 14 o^2 (-1 + Z^2) + 7 o^4 (-1 + Z^2) - o^6 (-1 + Z^2))=\eta_\infty(Mirror9_{42}).$$

\begin{figure}[h!]
\begin{center}
     \includegraphics[width=15cm]{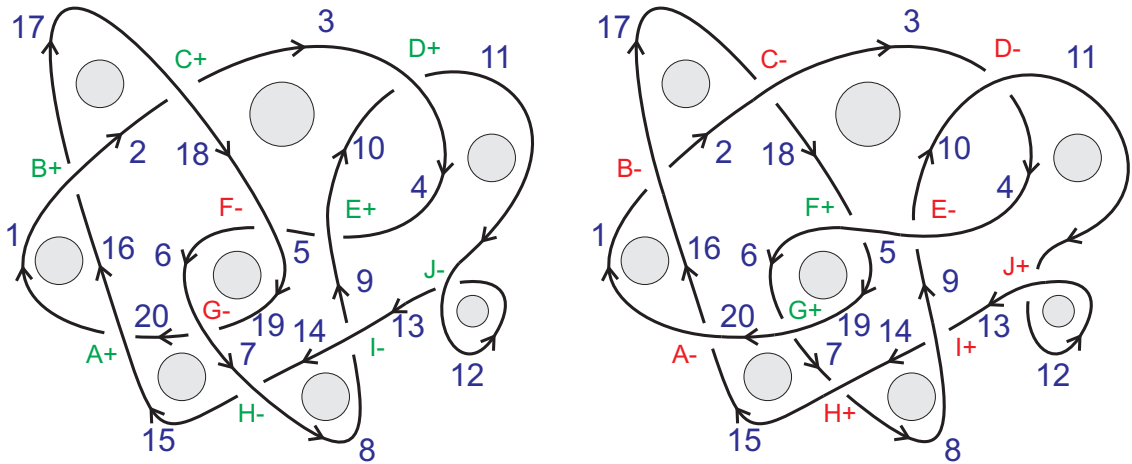}
         \caption{\bf Knots $9_{42}$ and $Mirror 9_{42}$: the 4 invariants agree, the knots are distinct}
     \label{fig:9_42andMirrorFromKnotilus}
\end{center}
\end{figure}

\subsection{Conway and Kinoshita-Terasaka knots and their doubled}
\begin{figure}[h!]
\begin{center}
     \includegraphics[width=10cm]{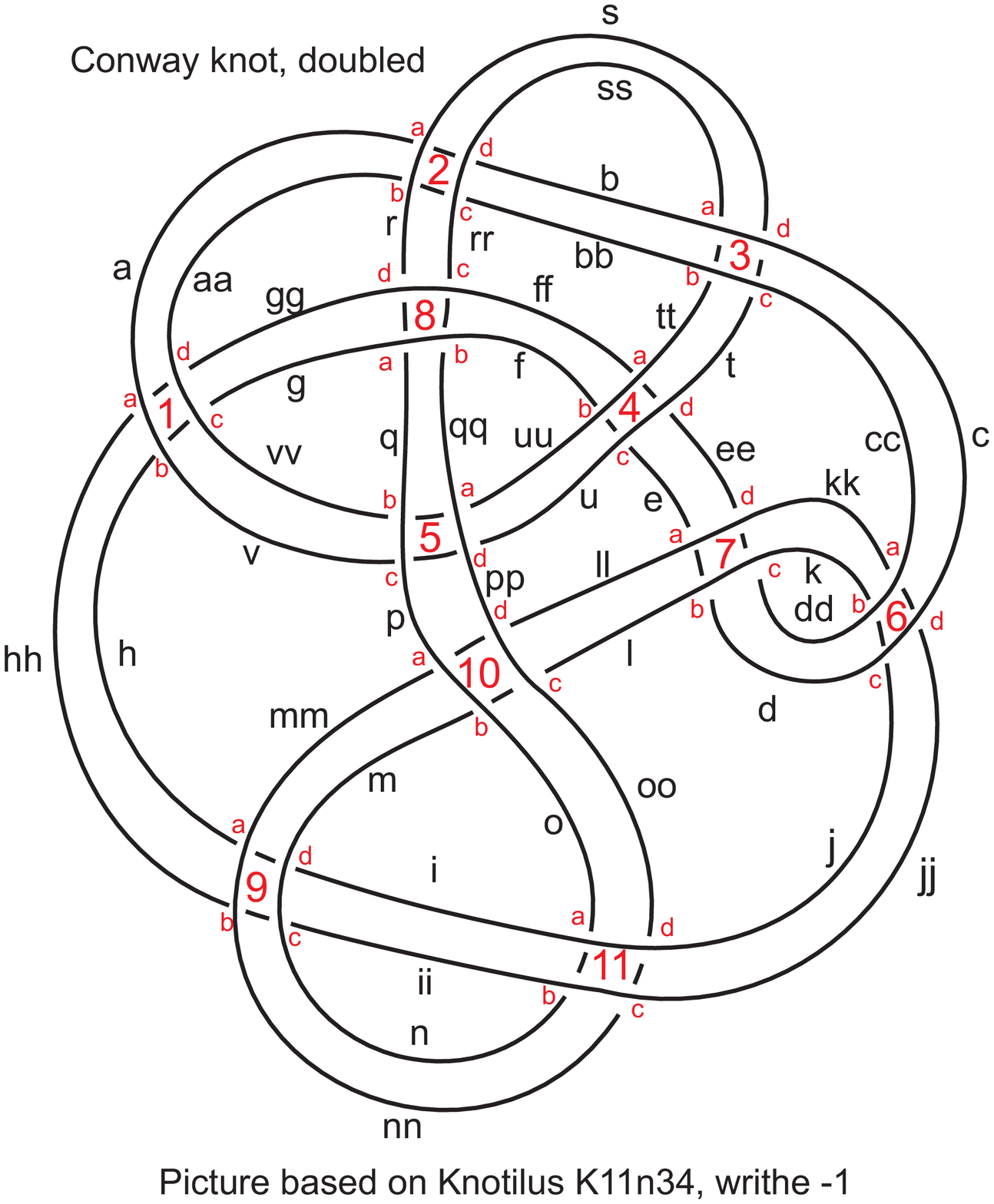}
         \caption{\bf Conway knot, doubled}
     \label{fig:ConwayDoubled}
\end{center}
\end{figure}

The few states of the specializations of the VSE-invariant for low level of $k$
can be used to effectively compute the doubled or tripled of some knots. I have done this for the Conway and Kinoshita-Terasaka knots hoping to detect mutation. The hope was not fulfilled because
for $k$ from 1 to 4 the VSE-specialiazation in these doubled knots are the same. Nevertheless
the hope persists for a higher level of cabling, by triplication or quadruplicating the knots. I show a summary of the computations. The links have each 44 crossing and is impossible to compute their bracket. For
the Conway doubled the computations start by encoding the crossings directly from the picture as follows:

\begin{verbatim}
c1 = X2[a, hh, ab1, ad1] X1[v, bc1, ab1, h] X2[vv, g, cd1, bc1]
   X1[aa, ad1, cd1, gg];
c2 = X2[s, a, ab2, ad2] X1[r, bc2, ab2, aa] X2[rr, bb, cd2, bc2]
   X1[ss, ad2, cd2, b];
c3 = X1[b, ab3, ad3, ss] X2[bb, tt, bc3, ab3] X1[cc, cd3, bc3, t]
   X2[c, s, ad3, cd3];
c4 = X2[tt, ff, ab4, ad4] X1[ab4, f, uu, bc4] X2[u, e, cd4, bc4]
   X1[t, ad4, cd4, ee];
c5 = X1[qq, ab5, ad5, uu] X2[q, vv, bc5, ab5] X1[p, cd5, bc5, v]
   X2[pp, u, ad5, cd5];
c6 = X2[cc, kk, ab6, ad6] X1[dd, bc6, ab6, k] X2[d, j, cd6, bc6]
   X1[c, ad6, cd6, jj];
c7 = X1[ll, ab7, ad7, e] X2[l, d, bc7, ab7] X1[k, cd7, bc7, dd]
   X2[kk, ee, ad7, cd7];
c8 = X2[g, q, ab8, ad8] X1[f, bc8, ab8, qq] X2[ff, rr, cd8, bc8]
   X1[gg, ad8, cd8, r];
c9 = X2[mm, h, ab9, ad9] X1[ab9, hh, nn, bc9] X2[n, ii, cd9, bc9]
   X1[m, ad9, cd9, i];
c10 = X2[p, mm, ab10, ad10] X1[ab10, m, o, bc10]
   X2[oo, l, cd10, bc10] X1[pp, ad10, cd10, ll];
c11 = X1[i, ab11, ad11, o] X2[ii, n, bc11, ab11]
   X1[jj, cd11, bc11, nn] X2[j, oo, ad11, cd11];
doubleConwayAsProduct = c1 c2 c3 c4 c5 c6 c7 c8 c9 c10 c11;
doubleConway = new[Link, doubleConwayAsProduct];
\end{verbatim}

Next, the normal forms are computed from the Mathematical command
\begin{center}
{\sf sumOfPolynomialNumberOfStates[aLink, k]]]}
\end{center} which I have implemented:

\begin{verbatim}
\[Eta]1doubleConway =
 Simplify[
  Subscript[\[Eta], 1][
   sumOfPolynomialNumberOfStates[doubleConway, 1]]]
\[Eta]2doubleConway =
 Simplify[
  Subscript[\[Eta], 2][
   sumOfPolynomialNumberOfStates[doubleConway, 2]]]
\[Eta]3doubleConway =
 Simplify[
  Subscript[\[Eta], 3][
   sumOfPolynomialNumberOfStates[doubleConway, 3]]]
\[Eta]4doubleConway =
 Simplify[
  Subscript[\[Eta], 4][
   sumOfPolynomialNumberOfStates[doubleConway, 4]]]
\end{verbatim}

Up to the level $k=4$ the values of the specializations of the VSE-invariant of Conway, doubled and Kinoshita-Terasaka, doubled are the same:

\noindent
$\eta_1(doubleConway)=o (-2 M o^2 Y Z + o (3 + 2 M Y Z - 2 Z^2) + 2 (-1 + Z^2))=\eta_1(doubleKT),$ \\
$\eta_2(doubleConway)=-(1/2) o (o^2 (5 + 4 M Y Z - 5 Z^2) - (-5 + 4 M Y Z - Z^2) (-1 +
      Z^2)+o (-2 + Z^2 - Z^4 + 4 M Y Z (-2 + Z^2)))=\eta_2(doubleKT),$\\
$\eta_3(doubleConway)=1/4 o ((-1 + Z^2) (-9 + 26 M Y Z - 4 Z^2 - 6 M Y Z^3 + Z^4) -
   2 o^2 (5 - 5 Z^2 + M Y (Z + 3 Z^3)) + o (5 - 5 Z^2 + 5 Z^4 - Z^6 + 2 M Y Z (14 - 13 Z^2 + 3 Z^4)))=\eta_3(doubleKT),$\\
$\eta_4(doubleConway)=-(1/32) o (4 o^2 (5 + 4 M Y Z - 5 Z^2) (1 + 3 Z^2) + (-1 +
      Z^2) (-53 + 167 Z^2 - 23 Z^4 + 5 Z^6 -
      8 M Y Z (31 - 16 Z^2 + 5 Z^4)) +
   o (-105 + 180 Z^2 - 130 Z^4 + 28 Z^6 - 5 Z^8 +
      8 M Y Z (-33 + 41 Z^2 - 21 Z^4 + 5 Z^6)))=\eta_4(doubleKT).$

\begin{figure}[h!]
\begin{center}
     \includegraphics[width=10cm]{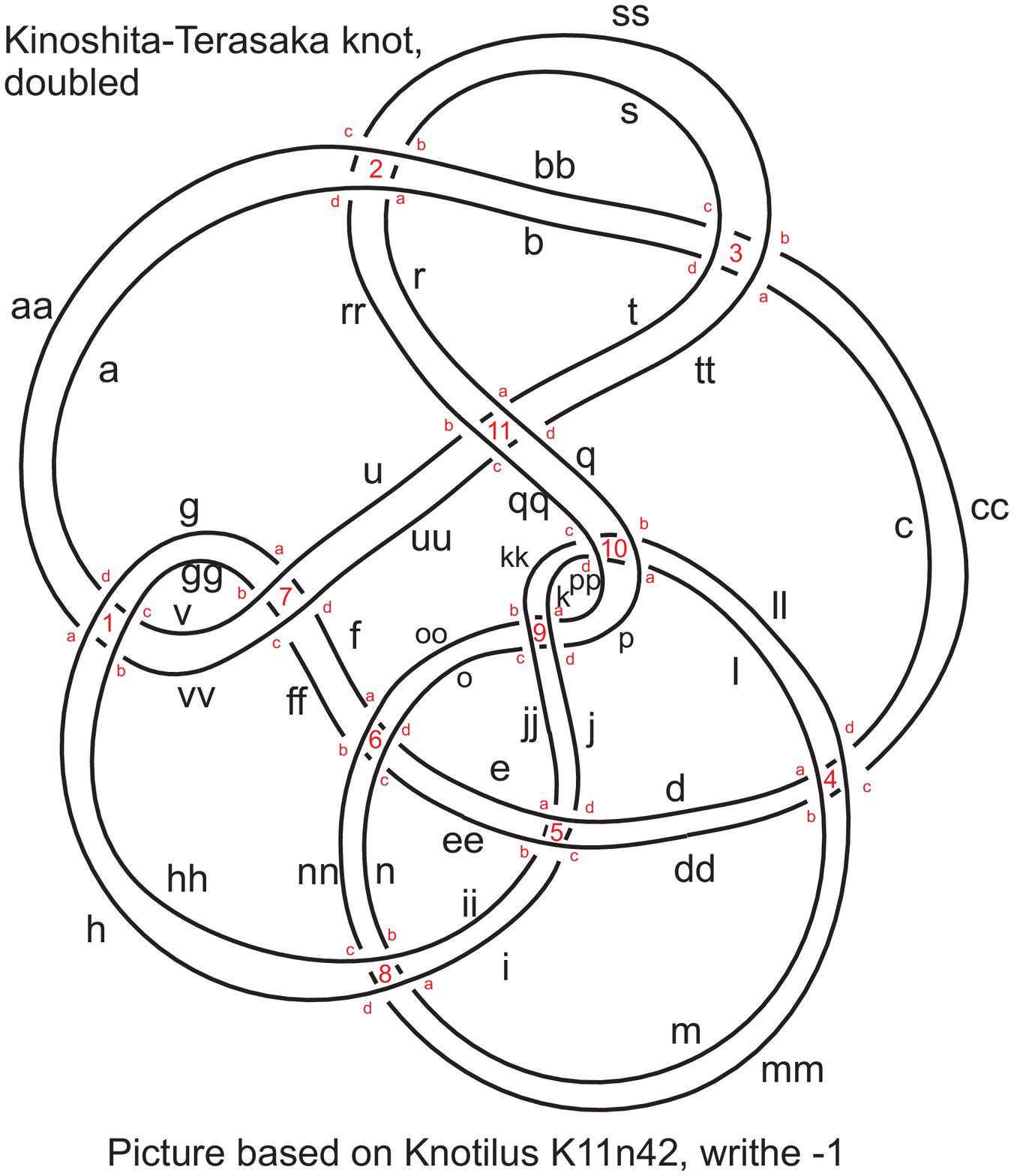}
         \caption{\bf Kinoshita-Terasaka knot, doubled}
     \label{fig:KTdoubled}
\end{center}
\end{figure}
\vspace{5mm}
\noindent
The computations for the doubled of Kinoshita-Terasaka knot follow similar lines.

\subsection{Links $T_{15}$ and $JS_{14}$}
\begin{figure}[h]
\begin{center}
     \includegraphics[width=12cm]{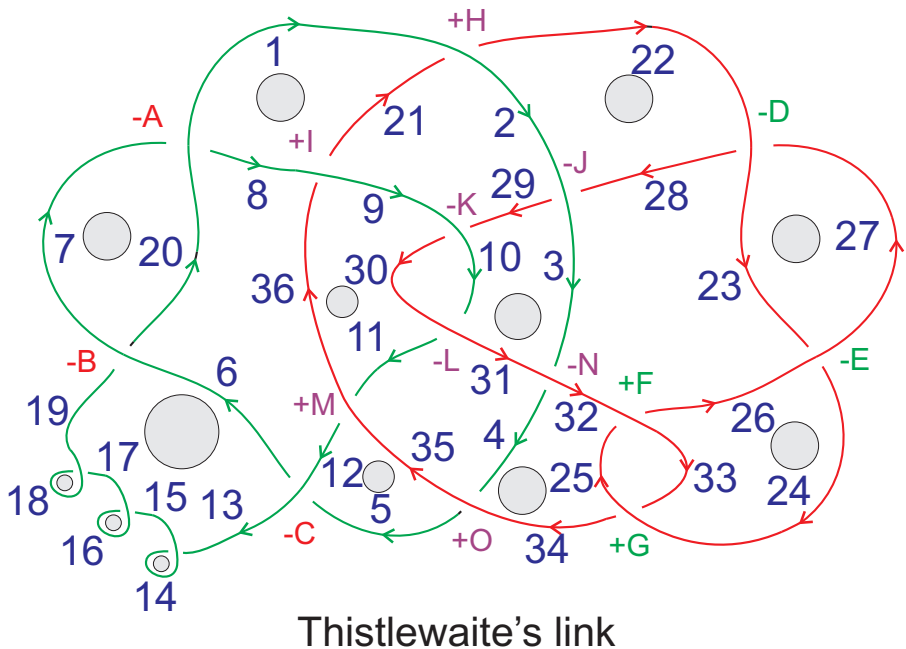}
         \caption{\bf Jones and VSE-invariant do not
         distinguish $T_{15}$ from the unlink} \label{fig:Thistle}
\end{center}
\end{figure}

Consider the link $T_{15}$, depicted in Fig.\ref{fig:Thistle} which is the first example of \cite{MBT2001}, writhe normalized to 0 at each component. The bracket polynomial of this link is equal to the bracket polynomial of the unlink: both are equal to o. Up to the level $k=4$, the VSE-invariant does not distinguishes $T_{15}$ from the unlink. The values of $\eta_i(T_{15})$ are all equal to $o^2$, for $i=1,2,3,4.$

\begin{figure}[h]
\begin{center}
     \includegraphics[width=12cm]{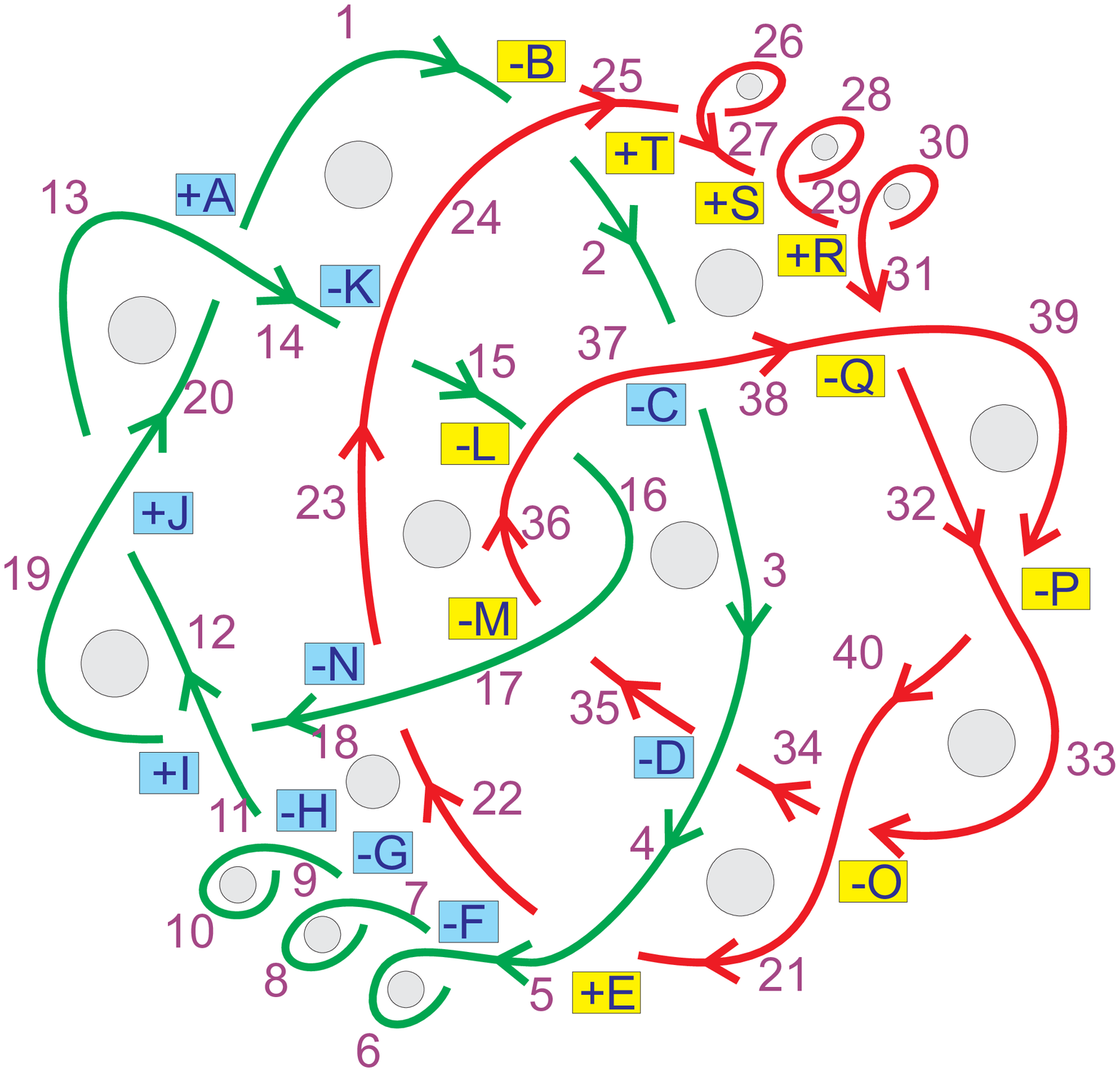}
         \caption{\bf Link $JS_{14}$ proves that the
         VSE-invariant is stronger than Jones'}\label{fig:JablanSazdanovic}
     \end{center}
\end{figure}

Consider also the link $JS_{14}$ in Fig. \ref{fig:JablanSazdanovic} based on a picture obtained from \cite{JS}. The bracket polynomial of this link is, once more, equal to the bracket polynomial of the unlink.

Here is a session of Mathematica computing, from Fig. \ref{fig:JablanSazdanovic}, the invariants $\eta_1, \eta_2, \eta_3, \eta_4$ for the link $JS_{14}$.
\begin{verbatim}
(* Crossing encoding for $JS_{14}$ *)
p1 = X1[13, 20, 14, 1] X2[25, 1, 24, 2] X1[38, 2, 37, 3]
     X1[3, 35, 4, 34] X2[4, 22, 5, 21];
p2 = X1[5, 7, 6, 6] X1[7, 9, 8, 8] X1[9, 11, 10, 10]
     X1[11, 18, 12, 19] X1[20, 13, 19, 12];
p3 = X1[24, 14, 23, 15] X2[37, 15, 36, 16] X2[16, 36, 17, 35]
     X1[17, 23, 18, 22] X2[40, 34, 21, 33];
p4 = X2[33, 39, 32, 40] X2[38, 32, 39, 31] X2[30, 29, 31, 30]
     X2[28, 27, 29, 28] X2[26, 25, 27, 26];
JSAsProduct = p1 p2 p3 p4;
JS = new[Link, JSAsProduct];
\eta1JS =
 Simplify[\eta1[sumOfPolynomialNumberOfStates[JS, 1]]]
\eta2JS =
 Simplify[\eta2[sumOfPolynomialNumberOfStates[JS, 2]]]
\eta3JS =
 Simplify[\eta3[sumOfPolynomialNumberOfStates[JS, 3]]]
\eta3JS =
 Simplify[\eta4[sumOfPolynomialNumberOfStates[JS, 4]]]
 \end{verbatim}
\noindent
Computations yield
$$\eta_1(JS) = o^2,$$
$$\eta_2(JS) = o \left(16 \left(Z^2-1\right) o^2+\left(16 Z^2-15\right) o-32
   \left(Z^2-1\right)\right),$$
$$\eta_3(JS) = o \left(16 \left(Z^2-1\right) o^2+\left(16 Z^2-15\right) o-32
   \left(Z^2-1\right)\right),$$
$$\eta_4(JS) = o (o^7+1024 Z^2 (1-2 Z^2)^2 (Z^2-1)
   o^2+1024 Z^2 (1-2 Z^2)^2 (Z^2-1) o- $$
$$ 2048
   Z^2 (1-2 Z^2)^2 (Z^2-1)).$$

Therefore $JS_{14}$ is distinguished from the unlink with two components and from $T_{15}$ at the $k=2$ specialization! This specialization has only 801 states as compared with the $2^{20}=1,048,576$ states of the Jones polynomial and the $3^{20}=3,486,784,401$ states of the full VSE-invariant.

\section{Conclusion and further research along these lines}
I have introduced a new regular isotopy invariant, the VSE-invariant which is a generalization of the Jones invariant (equivalent to Kauffman's bracket). The VSE-invariant is the normal form, $\eta_\infty$, of an class of polynomials $p \in \R/\I_\infty$, where $\R$ is the ring $$\R=\IZ[A,B,F,X,Y,Z,M,o]$$ and $\I_\infty$ is an ideal generated by 27 polynomials. The polynomial $p$ is computed directly from a diagram for the link: it is a state sum based on a virtual, shaded, exterior expansion (VSE-expansion) of the link diagram.

The scheme admits an infinite sequence of specializations based on the normal forms, $\eta_k$, k=1,2,\ldots. These are normal forms relative to a fixed Gr\"obner basis of the class of polynomials in $\R/\I_k$, where $\I_k = \langle \I_\infty \cup M^{k+1}\rangle$. Each $\eta_k$ induces a specialization of the VSE-invariant which is also a regular isotopy invariant.

I have proved that VSE-invariant is a strict generalization of the bracket.
For a fixed number of crossings $n$ the number of states to compute $\eta_k$ is a polynomial of degree $k$.  The specialization can be useful even in the case $k=2$, where the VSE-invariant distinguishes pairs of links not distinguishable by the bracket. At the low level $k=2,3$ the VSE-specialization can be computed obtaining useful invariants for links with hundred of crossings.

The possibility to detect mutants via $m$-cabling failed for the Conway and Kinoshita-Terasaka knots at the level $m=2$. But it might work for higher values of $m$. Testing this idea awaits for a proper implementation of $m$-cabling.

Variations of the strategy here introduced can generalize the VSE-invariant and be exported to obtain invariants of 3-manifolds with a polynomial number of states along the combinatorial lines of \cite{KP} and \cite{KL}. These matters will be treated in future papers.

\section{Acknowledgement} My research is supported by a Grant of CNPq, Brazil (process number 306106/2006) and by a University Position at UFPE, Recife, Brazil.

\end{document}